\documentclass[10pt,a4paper]{article}

\usepackage[utf8]{inputenc}
\usepackage[english]{babel}

\usepackage[T1]{fontenc}
\usepackage[english]{babel}

\usepackage{amsmath}
\usepackage{amsfonts}
\usepackage{amssymb}
\usepackage{amsthm}

\usepackage{graphicx}
\usepackage[margin=2cm]{geometry}
\usepackage{paralist}

\usepackage{csquotes,etoolbox,keyval,ifthen,url,babel}
\usepackage[backend=biber, defernumbers=true, sorting=anyt, style=numeric-comp, maxbibnames=4, maxcitenames=2]{biblatex}

\renewbibmacro{in:}{%
  \ifentrytype{article}{}{\printtext{\bibstring{in}\intitlepunct}}}
\AtBeginDocument{
                 }

\usepackage{url}                         
\usepackage[breaklinks]{hyperref}        %
\usepackage{breakurl}                    %
\hypersetup{                                
    colorlinks=true,         
    linkcolor=red,
    filecolor=cyan,      
    urlcolor=cyan,
    breaklinks=true
            } 

\usepackage{expl3}                                                                       %
\ExplSyntaxOn                                                                            %
\newcommand \breakDOI[1]                                                                 %
   {                                                                                     %
	 \tl_map_inline:nn { #1 } { \href{http://dx.doi.org/#1}{##1} \penalty0 \scan_stop: } %
   }                                                                                     %
\ExplSyntaxOff                                                                           %
\DeclareFieldFormat{doi}{
	\mkbibacro{DOI}\addcolon\space                                                       %
	{\breakDOI{#1}}}                                                                     %
\usepackage{amsthm}
\usepackage{chngcntr}

\usepackage{caption}
\usepackage{subcaption}
\usepackage{wrapfig}

\usepackage[font=small,labelfont=bf,tableposition=top]{caption}
\DeclareCaptionLabelFormat{andtable}{#1~#2  \&  \tablename~\thetable}

\usepackage{booktabs}

\newtheorem{theorem}{Theorem}
\newtheorem{lemma}{Lemma}

\newtheorem{example}{Example}
\newtheorem{remark}{Remark}

\counterwithin{theorem}{section}
\counterwithin{definition}{section}
\counterwithin{example}{section}
\counterwithin{lemma}{section}
\counterwithin{remark}{section}

\DeclareMathOperator{\dif}{d}

\graphicspath{{Pictures/}}

\allowdisplaybreaks

\title{A fitted second--order difference scheme on a modified Shishkin mesh for a semilinear singularly-perturbed boundary-value problem}
\author{Samir Karasulji\'c\footnote{corresponding author}\:$^{,\thinspace}$\footnote{University of Tuzla, Faculty of Sciences and Mathematics, 
		Univerzitetska br. 4, 75 000 Tuzla, Bosnia and Herzegovina, email:\texttt{samir.karasuljic@unitz.ba}} 
	and Irma Zenunovi\'c \footnote{College of finance and accounting FINra Tuzla, Mitra Trifunovi\'ca U\v ce 9, 75 000 Tuzla, Bosnia and Herzegovina, email: \texttt{irma.finra@gmail.com}}}
\date{}

\begin{document}

\maketitle 
	
\begin{abstract}
	In the present paper we consider the numerical solving of a semilinear singular--perturbation reaction--diffusion boundary--value problem having boundary layers. A new difference scheme is constructed, the second order of convergence on a modified Shishkin mesh is shown.  The numerical experiments are included in the paper,  which confirm the theoretical results. 
\end{abstract}
	
\section{Introduction}
Let us consider the following boundary value problem
\begin{gather}
  \varepsilon^2y''=f(x,y),\quad x\in(0,1),\label{uvod1}\\
   y(0)=y(1)=0,\label{uvod2}
\end{gather}
where $0<\varepsilon<1,$ and we assume that the nonlinear function $f$ is continuously differentiable, i.e. for $k\geqslant 2,$ $f\in C^{k}([0,1]\times\mathbb{R}),$ and that it has a strictly positive derivative with respect  to $y$ 
\begin{equation}\label{uvod3}
 \frac{\partial f}{\partial y}:= f_y\geqslant m>0,\quad (x,y)\in[0,1]\times\mathbb{R},
\end{equation}
where $m$ is a constant. The problem \eqref{uvod1}--\eqref{uvod2} under the condition \eqref{uvod3} has the unique solution, see Lorentz \cite{lorenz1982stability}. It's a well known fact that the exact solution to the problem \eqref{uvod1}--\eqref{uvod3} rapidly changes near the ends points $x=0$ and $x=1.$ Also, when any classical numerical method is applied to the problem \eqref{uvod1}--\eqref{uvod3},  then a large number of mesh points must be used to obtain a satisfactory numerical solution. It is very expensive from the computing side, and that is one of reasons for developing numerical methods which taking into account the properties of singular--perturbation boundary--value problems.

Many authors have worked on the numerical treatment of the problem \eqref{uvod1}--\eqref{uvod3} with different assumptions about the function $f,$  and as well as more general nonlinear problems. There were many constructed $\varepsilon$--uniformly convergent difference schemes of order 2 and higher (Herceg \cite{herceg1990}, Herceg, Surla and Rapaji{\'c} \cite{herceg1991}, Herceg and Miloradovi{\'c} \cite{herceg2003}, Herceg and Herceg \cite{herceg2003a}, Kopteva and Lin\ss\, \cite{kopteva2001}, Kopteva and Stynes \cite{kopteva2001robust,kopteva2004}, Kopteva, Pickett and Purtill \cite{kopteva2009}, Lin\ss, Roos and Vulanovi{\'c} \cite{linss2000uniform}, Sun and Stynes \cite{stynes1996}, Stynes and Kopteva \cite{stynes2006numerical}, Surla and Uzelac \cite{surla2003}, Vulanovi{\'c} \cite{vulanovic1983, vulanovic1989, vulanovic1991second, vulanovic1993, vulanovic2004}, Kopteva \cite{kopteva2007maximum} etc.). 

In this paper we use a method introduced by Boglaev  \cite{boglaev1984approximate} to construct a different scheme, and a layer--adapted mesh introduced by Vulanovi\' c  \cite{vulanovic2001higher}. In the paper \cite{samir2015uniformly} the difference scheme was constructed for the same problem and the authors used the same layer--adapted mesh. It is well known fact that Shishkin mesh is the simplest layer--adapted mesh. It looks like as two uniform meshes glued, the one finer the other coarser. The simplicity of this mesh implies  some simpler analysis than using other layer--adapted meshes, and that is the main benefit of using Shishkin mesh.  Unfortunately, the cost of simplicity is a greater value of error.

The paper is organized as follows. The first section is Introduction, where the problem is listed and the main results. Layer--adapted mesh is 2nd section, here are given generating functions of meshes we use in this paper. A new difference scheme is constructed in section The difference scheme. Their stability is shown in the section Stability.  The uniform convergence is proven in 6th section, 7th section is Numerical experiments and the last section is Conclusion.

The author's results in numerical solving the problem \eqref{uvod1}--\eqref{uvod3} and others results can be seen in 
\cite{samir2010scheme}, \cite{samir2011scheme}, \cite{samir2015uniformlyconvergent}, \cite{samir20152d}, \cite{samir2015uniformly}, \cite{samir2015modifiedbah}, \cite{samir2017construction}, \cite{samir2018uniformly}, \cite{karasuljic2019class}, \cite{liseikin2020numerical}, \cite{karasuljic2020onconstruction}, \cite{Liseikin_2021}, \cite{samir2010tmt}, \cite{samir2011uniformnly}, \cite{samir2011skoplje}, \cite{samir2012uniformnly}, \cite{samir2012class}, \cite{samir2013collocation}, \cite{samir2015construction}, \cite{liseikin2019rules}, \cite{liseikin2020comprehensiveproceedings}, \cite{samir2021matematika1}, \cite{liseikin2021numerical}.

\section{Layer--adapted meshes}\label{mesh}

We mentioned that classical methods   are not suitable for problems like \eqref{uvod1}--\eqref{uvod2}, before. There are several issues, with the stability, the cost of calculations and so on. These issues have been overcome by constructing special methods which take into account the presence of a layer or layers. One such a method is the fitted mesh method. Meshes that were constructed by this method have the nonuniform distribution of mesh points. The distribution of points is dictated by the behavior of the exact solution and their derivatives in a layer or layers. Estimates given in the following theorem are key to understanding of this behavior. In the following analysis we need the decomposition of the solution $y$ of the problem 
$(\ref{uvod1})-(\ref{uvod2})$ to a layer component $s$ and a regular component $r$, which is given in the following theorem.

\begin{theorem} {\rm\cite{vulanovic1983} \label{theoremDecomposition}
	} The solution $y$ to the problem $(\ref{uvod1})-(\ref{uvod2})$ can be represented in the following way:
	\begin{equation*}
	y=r+s,
	\end{equation*}
	where for $j=0,1,...,k+2$ and $x\in[0,1]$ we have  that
	\begin{equation}
	\left|r^{(j)}(x)\right|\leq  C,
	\label{regularna}
	\end{equation}
	and
	\begin{equation}
	\left|s^{(j)}(x)\right|\leq  C \varepsilon^{-j}\left(e^{-\frac{x}{\varepsilon}\sqrt{m}}+e^{-\frac{1-x}{\varepsilon}\sqrt{m}}\right).
	\label{slojna}
	\end{equation}
	\label{teorema1}
\end{theorem}

Based on the previous theorem, it's a well--known fact that the exact solution to problem \eqref{uvod1}--\eqref{uvod3} changes rapidly near the end points $x=0$ and $x=1.$

Many meshes have been constructed for the numerical solving problems having a layer or layers of an exponential type. In the present paper we shall use four different meshes. Let $N+1$ be the number of mesh points. These meshes $0=x_0<x_1<\ldots<x_N=1,$ we will get by using appropriate generating functions, i.e. $x_i=\psi(i/N).$ The generating functions are constructed as follows. 

The first mesh is Shishkin mesh \cite{shishkin1988grid}, the generating function for this mesh is 
\begin{equation}\label{meshShishkin1}
\psi(t)=\begin{cases}
          4\lambda  t,\quad t\in[0,1/4] \\
           \lambda+2(1-2\lambda)(t-1/4),\quad t\in[1/2,1/4],\\
           1-\psi(1-t),\quad t\in[1/2,1],
        \end{cases}
\end{equation}
where $\lambda$ Shishkin mesh transition point by
\begin{equation}\label{shihskinTransitionPoint}
\lambda:=\min\left\{\frac{2\varepsilon\ln N}{\sqrt{m}},\frac{1}{4}\right\} .
\end{equation}

The second mesh is modified Shishkin mesh proposed by Vulanovi\'c \cite{vulanovic2001higher}, the generating function for this mesh is 
\begin{equation}\label{meshShishkin2}
\psi(t)=\begin{cases}
          4\lambda t,\quad  t\in[0,1/4],\\
          p(t-1/4)^3+4\lambda t, \quad t\in[1/4,1/2],\\
          1-\psi(1-t),\quad t\in[1/2,1],
        \end{cases}  
\end{equation}
where $p$ is chosen so that $\psi(1/2)=1/2,$ i.e. $p=32 (1-4\lambda).$ Note that $\psi\in C^{1}[0,1]$ with $\Vert\psi'\Vert_{\infty}\leqslant C,$ $\Vert\psi''\Vert_{\infty}\leqslant C.$ Therefore the mesh size $h_i=x_{i+1}-x_i,\,i=0,\ldots N-1$ satisfy (see \cite{linss2012approximation})
\begin{eqnarray}
h_i=\int_{i/N}^{(i+1)N}{\psi'(t)\dif t}\leqslant CN^{-1},\quad |h_{i+1}-h_i|=\left|\int_{i/N}^{(i+1)/N}{\int^{t+1/N}_{t}{\psi''(s) \dif s}}\right|\leqslant CN^{-2}.
\end{eqnarray}
The Shishkin mesh transition point $\lambda$ is the same as in the first Shishkin mesh, i.e. \eqref{shihskinTransitionPoint}.

The third mesh is modified Bakhvalov mesh also proposed by Vulanovi\'c \cite{vulanovic1983numerical}, the generating function for this mesh is 

\begin{equation}\label{meshVulanovic1}
     \psi(t)=\begin{cases}
               \mu(t):=\frac{a\varepsilon t}{q-t},\quad t\in[0,\alpha],\\
               \mu(\alpha)+\mu'(\alpha)(t-\alpha),\quad t\in[\alpha ,1/2],\\
               1-\psi(1-t),\quad t\in[1/2,1],
             \end{cases}
\end{equation}
where $a$ and $q$ are constants, independent of $\varepsilon,$ such that $q\in(0,1/2),\:a\in(0,q/\varepsilon),$ and additionally $a\sqrt{m}\geqslant 2.$ The parameter 
$\alpha$ is the abscissa of the contact point of tangent line from $(1/2,1/2)$ to $\mu (t),$ and its value is 
\[\alpha =\frac{q-\sqrt{aq\varepsilon(1-2q+2a\varepsilon)}}{1+2a\varepsilon }.\]

The last but not least is the mesh proposed by Liseikin \cite{liseikin2018grid, liseikin2019compact}, and we will use 	its modification from \cite{liseikin2020numerical}. The generating function for this mesh is 

\begin{equation}\label{meshLiseikin1}
\psi( t,\varepsilon,a,k)=
\left\{
\begin{array}{ll}
\displaystyle c_1\varepsilon^k ((1-d t) ^{-1/a}-1)\;, & 0\leqslant t \leqslant 1/4 \;,\\[4mm]
\displaystyle c_1\Bigl[\varepsilon^{kan/(1+na)}-\varepsilon^k+d\frac{1}{a}\varepsilon^{ka(n-1)/(1+na)}(t-1/4 )+\\
\frac{1}{2}d^2\frac{1}{a}\Bigl(\frac{1}{a}+1\Bigr)\varepsilon^{ka(n-2)/(1+na)}(t -1/4 )^2+c_0(t -1/4 )^{3}\Bigr]\;, & 1/4  \leqslant  t \leqslant 1/2\;,\\
1-\psi(1-t ,\varepsilon,a,k)\;,& 1/2\leqslant t \leqslant 1\;,
\end{array}
\right .
\end{equation}
where $d=(1-\varepsilon^{ka/(1+na)})/(1/4),$ $a$ is a positive constant subject to $a\geq m_1 > 0$, and $a=1,$  $c_0>0$, $n=2,$ $k=1,$ $c_0=0,$   and
$
\tfrac{1}{c_1}=2\left[\varepsilon^{kan/(1+na)}-\varepsilon^k+\tfrac{d}{4a} \varepsilon^{ka(n-1)/(1+na)}\right.
\left.+ \tfrac{d^2}{2}\tfrac{1}{a}\Bigl(\frac{1}{a}+1\Bigr)\varepsilon^{ka(n-2)/(1+na)}(1/4 )^2+c_0(1/4)^{3}\right]
$
is chosen here.

\section{The difference scheme}
The first step in the numerical solving of the problem \eqref{uvod1}--\eqref{uvod3} is a construction of  difference scheme, which generates a nonlinear system of equations. A solution of this nonlinear system  is a discrete numerical solution to the problem \eqref{uvod1}--\eqref{uvod3}. 
\subsection{Construction of the difference scheme}

From the paper \cite{boglaev1984approximate} we have the following equality   
\begin{multline}
\frac{\beta}{\sinh(\beta h_{i-1})}y_{i-1}-\left(\frac{\beta}{\tanh(\beta h_{i-1})}+\frac{\beta}{\tanh(\beta h_i)} \right)y_i+\frac{\beta}{\sinh(\beta h_i)}y_{i+1}=
       \\\frac{1}{\varepsilon^2}\left[\int\limits_{x_{i-1}}^{x_{i}}{u_{i-1}^{II}(s)\psi(s,y(s))\dif s}+\int\limits_{x_{i}}^{x_{i+1}}{u_{i}^{I}(s)\psi(s,y(s))\dif s}\right],
\label{konst1}
\end{multline}
\[y_0=0,\:\:y_N=0,\:i=1,2,\ldots,N-1,\]
where 
\[u^{I}_{i}(x)=\frac{\sinh(\beta (x_{i+1}-x))}{\sinh(\beta h_{i})},\,u^{II}_i(x)=\frac{\sinh(\beta (x-x_i))}{\sinh(\beta h_i)},\:x\in[x_i,x_{i+1}].\]
In the general case, we cannot explicitly calculate integrals on the right hand side \eqref{konst1}.  In dependence how we approximate the integrals in \eqref{konst1} we get various difference schemes. In the papers \cite{samir2015uniformlyconvergent, samir2015uniformly, karasuljic2019class}, the approximates are not so simple. In the papers \cite{samir2015uniformlyconvergent, samir2015uniformly, karasuljic2019class}, the approximations of integrals could be  simpler, therefore the analysis of methods would be easier. In order to avoid any difficulties  and make the analysis easier as we can , we shall use  the following approximation for the integrals in \eqref{konst1}:   we approximate the function $\psi$ on the intervals $[x_{i-1},x_i],$ $i=1,2,\ldots, N$ by 
\begin{equation}
 \overline{\psi}_i=\frac{\psi(x_{i-1},\overline{y}_{i-1})+\psi(x_i,\overline{y}_i)}{2},
\label{konst2}
\end{equation} 
where $\overline{y}_i$ is an approximate value of the solution of the problem \eqref{uvod1}--\eqref{uvod3} at points $x_i,$ and $\psi(x,y)=f(x,y)-\gamma y$ $\gamma$ is such a constant, that  holds
\begin{equation}
\gamma\geqslant  f_y.
\label{konst3}
\end{equation}

Putting \eqref{konst2} into \eqref{konst1}, and after some computation, we get the following difference scheme

\begin{multline}  
\frac{\cosh(\beta h_{i-1})+1}{\sinh(\beta h_{i-1})}\overline{y}_{i-1}-\left(\frac{\cosh(\beta h_{i-1})+1}{\sinh(\beta h_{i-1})} 
     +\frac{\cosh(\beta h_{i})+1}{\sinh(\beta h_{i})}\right)\overline{y}_i+\frac{\cosh(\beta h_{i})+1}{\sinh(\beta h_{i})}\overline{y}_{i+1}\\
     -\frac{f_{i-1}+f_i}{\gamma}\cdot\frac{\cosh(\beta h_{i-1})-1}{\sinh(\beta h_{i-1})}-\frac{f_{i}+f_{i+1}}{\gamma}\cdot\frac{\cosh(\beta h_{i})-1}{\sinh(\beta h_{i})}=0,
\label{konst4}     
\end{multline}
where  $f_i=f(x_i,\overline{y}_i),$ $\overline{y}_0=y(0),\overline{y}_N=y(1),$ and $\:i=1,\ldots,N-1.$ 

\section{Stability}\label{stability}
The difference scheme \eqref{konst4} generates a system of nonlinear algebraic equations. A solution of this system is a discrete numerical solution of the problem \eqref{uniform3}--\eqref{uvod3}.  The next tasks are to show the existence and uniqueness of the discrete numerical solution and the stability of the difference scheme. Let us set the discrete operator 
\begin{equation}
    Tu=(T u_0,Tu_1,\ldots,Tu_N)^T,
 \label{operator1}
\end{equation}
where 
\begin{align}
    Tu_0&=-u_0\nonumber\\
    Tu_i&=\frac{\gamma}{\frac{\cosh(\beta h_{i-1})-1}{\sinh(\beta h_{i-1})}+\frac{\cosh(\beta h_{i})-1}{\sinh(\beta h_{i})}}  \nonumber\\
    &\cdot\left[\frac{\cosh(\beta h_{i-1})+1}{\sinh(\beta h_{i-1})}u_{i-1}-\left(\frac{\cosh(\beta h_{i-1})+1}{\sinh(\beta h_{i-1})} 
     +\frac{\cosh(\beta h_{i})+1}{\sinh(\beta h_{i})}\right)u_i+\frac{\cosh(\beta h_{i})+1}{\sinh(\beta h_{i})}u_{i+1}\right.\nonumber\\
     &\quad-\left.\frac{f(x_{i-1},u_{i-1})+f(x_i,u_i)}{\gamma}\cdot\frac{\cosh(\beta h_{i-1})-1}{\sinh(\beta h_{i-1})}-\frac{f(x_{i},u_i)
             +f(x_{i+1},u_{i+1})}{\gamma}\cdot\frac{\cosh(\beta h_{i})-1}{\sinh(\beta h_{i})}\right],\\
     &\hspace{5cm}    i=1,\ldots,N-1\nonumber \\
    Tu_N&=-u_N. \nonumber
\end{align}
Obviously, it is hold
\begin{equation}
    T\overline{y}=0,
 \label{operator2}
\end{equation}
where $\overline{y}=(\overline{y}_0,\overline{y}_1,\ldots,\overline{y}_N)^T$ the numerical solution of the problem \eqref{uvod1}--\eqref{uvod3}, obtained by using the difference scheme \eqref{konst4}.

\begin{theorem}\label{thereomStability1}
   The discrete problem \eqref{operator1}--\eqref{operator2} has a unique solution $\overline{y}$ for $\gamma\geqslant f_y.$ Moreover, for every $v,w\in\mathbb{R}^{N+1}$ we have the following stability inequality \[\Vert v-w\Vert\leqslant C\Vert Tv-Tw\Vert.\] 
\end{theorem}

\begin{proof}
Denote the Fr\' echet derivative of the discrete operator $T$ by $H,$ e.g $H:=(T\overline{y})',$ and $H=(h_{ij}).$ Now, the non-zeros  elements of this matrix $H$ are
    \begin{align*}
     h_{1,1}=&-1,\:h_{N+1,N+1}=-1,\\
     h_{i,i-1}=&\frac{\gamma}{\frac{\cosh(\beta h_{i-1})-1}{\sinh(\beta h_{i-1})}+\frac{\cosh(\beta h_{i})-1}{\sinh(\beta h_{i})}}
            \left( \frac{\cosh(\beta h_{i-1})+1}{\sinh(\beta h_{i-1})}
                    -\frac{\frac{\partial f}{\partial \overline{y}_{i-1}}}{\gamma}\frac{\cosh(\beta h_{i-1})-1}{\sinh(\beta h_{i-1})}\right),\\
     h_{i,i+1}=&\frac{\gamma}{\frac{\cosh(\beta h_{i-1})-1}{\sinh(\beta h_{i-1})}+\frac{\cosh(\beta h_{i})-1}{\sinh(\beta h_{i})}} 
                \left( \frac{\cosh(\beta h_{i})+1}{\sinh(\beta h_{i})}
                        -\frac{\frac{\partial f}{\partial \overline{y}_{i+1}}}{\gamma}\cdot \frac{\cosh(\beta h_{i})-1}{\sinh(\beta h_{i})}\right),\\
     h_{i,i}=&-\frac{\gamma}{\frac{\cosh(\beta h_{i-1})-1}{\sinh(\beta h_{i-1})}+\frac{\cosh(\beta h_{i})-1}{\sinh(\beta h_{i})}} 
              \left[\frac{\cosh(\beta h_{i-1})+1}{\sinh(\beta h_{i-1})}+\frac{\cosh(\beta h_{i})+1}{\sinh(\beta h_{i})}\right.\\       
             &\hspace{2cm}\left. +\frac{\cosh(\beta h_{i-1})-1}{\sinh(\beta h_{i-1})}\cdot\frac{\frac{\partial f}{\partial \overline{y}_{i}}}{\gamma} 
                       +\frac{\cosh(\beta h_{i})-1}{\sinh(\beta h_{i})} \cdot \frac{\frac{\partial f}{\partial \overline{y}_{i}}}{\gamma}\right],
                   \: i=2,\ldots,N.
    \end{align*}
Because \eqref{konst3},  we have 
\begin{equation}
   h_{1,1}<0,\:h_{N+1,N+1}<0,\:h_{i-1,i}\geqslant 0,\:h_{i+1,i}\geqslant 0,\:h_{i,i}<0,
\label{operator3}
\end{equation}
and 
\begin{equation}
   \left|h_{i,i}\right| -\left|h_{i,i-1}\right|-\left|h_{i,i+1}\right|\geqslant 2m>0,
 \label{operator4}
\end{equation}  
so we conclude that $H$ is an $M$--matrix, and 
\begin{equation}
    \Vert H^{-1}\Vert\leqslant C.
 \label{operator5}
\end{equation} 
Now, by Hadamard theorem \cite[Th 5.3.10]{ortega2000}, the discrete operator is  a homeomorphism, and \eqref{operator2} has the unique solution. 

The second statement of the theorem follows from
\[Tv-Tw=(T\xi)'(v-w), \]
for some $\xi=(\xi_0,\xi_1,\ldots,\xi_n)\in\mathbb{R}^{N+1},$ and based on \eqref{operator5} we finally get
\[\Vert v-w\Vert\leqslant C\Vert Tv-Tw\Vert.\]
\end{proof}
\begin{remark}
  If we the difference scheme \eqref{konst4} multiply by $-1,$ we will get 
   $h_{1,1}<0,\:h_{N+1,N+1}<0,\:h_{i,i-1}\geqslant 0,\:h_{i,i+1}\geqslant 0,\:h_{i,i}<0.$
\end{remark}

\section{Uniform convergence}

In this section we deal with a very important issue in the numerical solving of the problem \eqref{uvod1}--\eqref{uvod3}, that is error.  To proof Theorem on convergence we need the following lemmas.

\begin{lemma}\label{lema1}
  Assume that $\varepsilon \leqslant \frac{C}{N}.$ In the part of the modified mesh \eqref{meshShishkin2} from Section \ref{stability} when $x_i,\,x_{i\pm 1}\in[x_{N/4-1},\lambda]\cup[\lambda,1/2],$ we have the following estimate
 
\begin{multline*}\frac{\gamma}{\frac{\cosh(\beta h_{i-1})-1}{\sinh(\beta h_{i-1})}+\frac{\cosh(\beta h_i)-1}{\sinh(\beta h_i)}}
 \left|\frac{f(x_{i-1},y(x_{i-1}))+f(x_i,y(x_i))}{\gamma}\cdot \frac{\cosh(\beta h_{i-1})-1}{\sinh(\beta  h_{i-1})}\right.\\
              +\left.\frac{f(x_{i},y(x_{i}))+f(x_{i+1},y(x_{i+1}))}{\gamma}\cdot\frac{\cosh(\beta h_{i})-1}{\sinh(\beta  h_{i})}\right|
              \leqslant\frac{C}{N^2}, \:i=N/4,\ldots,N/2-1.
\end{multline*}    
\end{lemma}
\begin{proof}
   Due to $\varepsilon^2y''(x_i)=f(x_i,y(x_i)),$ theorem of decomposition for both components $r$ and $s$ in part of the mesh corresponding to  $[x_{N/4-1},\lambda]\cup[\lambda,1/2]$  and assumption $\varepsilon\leqslant \frac{C}{N},$ we have that
\begin{multline*}
\frac{\gamma}{\frac{\cosh(\beta h_{i-1})-1}{\sinh(\beta h_{i-1})}+\frac{\cosh(\beta h_i)-1}{\sinh(\beta h_i)}}\\
 \cdot\left| \frac{f(x_{i-1},y(x_{i-1}))+f(x_i,y(x_i))}{\gamma}\cdot \frac{\cosh(\beta h_{i-1})-1}{\sinh(\beta  h_{i-1})}
              +\frac{f(x_{i},y(x_{i}))+f(x_{i+1},y(x_{i+1}))}{\gamma}\cdot\frac{\cosh(\beta h_{i})-1}{\sinh(\beta  h_{i})}\right|\\
   \leqslant\left|f(x_{i-1},y(x_{i-1}))+f(x_i,y(x_i))\right|
                        +\left|f(x_{i},y(x_{i}))+f(x_{i+1},y(x_{i+1}))\right|\leqslant\frac{C}{N^2}.
\end{multline*} 
\end{proof}

\begin{lemma}{\rm\cite{samir2015uniformly}}\label{lema2}
  Assume that $\varepsilon \leqslant \frac{C}{N}.$ In the part of the modified mesh \eqref{meshShishkin2} from Section \ref{mesh} when $x_i,\,x_{i\pm 1}\in[x_{N/4-1},\lambda]\cup[\lambda,1/2],$ we have the following estimate
\begin{equation*}
\left| \frac{\frac{\cosh(\beta h_{i-1})-1}{\sinh(\beta h_{i-1})}(y(x_{i-1})-y(x_i))
              -\frac{\cosh(\beta h_{i})-1}{\sinh(\beta h_{i})}(y(x_{i})-y(x_{i+1})) }
              {\frac{\cosh(\beta h_{i-1})-1}{\sinh(\beta h_{i-1})}+\frac{\cosh(\beta h_{i})-1}{\sinh(\beta h_{i})}}\right|\leqslant\frac{C}{N^2},\:i=N/4,\ldots,N/2-1.
\end{equation*}
\end{lemma}

\begin{lemma}{\rm\cite{samir2015uniformly}}\label{lema3}
  Assume that $\varepsilon \leqslant \frac{C}{N}.$ In the part of the modified mesh \eqref{meshShishkin2} from Section \ref{mesh} when $x_i,\,x_{i\pm 1}\in[x_{N/4-1},\lambda]\cup[\lambda,1/2],$ we have the following estimate
\begin{equation*}
\left| \frac{\frac{y(x_{i-1})-y(x_i)}{\sinh(\beta h_{i-1})}
              -\frac{y(x_{i})-y(x_{i+1})}{\sinh(\beta h_{i})} }
              {\frac{\cosh(\beta h_{i-1})-1}{\sinh(\beta h_{i-1})}+\frac{\cosh(\beta h_{i})-1}{\sinh(\beta h_{i})}}\right|\leqslant\frac{C}{N^2},\:i=N/4,\ldots,N/2-1.
\end{equation*}
\end{lemma}
\noindent Now, we can state and prove the theorem on convergence. 

\begin{theorem}
The discrete problem \eqref{operator1}--\eqref{operator2} on the modified Shishkin mesh \eqref{meshShishkin2} from Section \ref{mreze} is uniformly convergent with respect $\varepsilon$ and 
\begin{equation*}
  \max_i\left|y(x_i)-\overline{y}_i\right|\leqslant C 
         \left\{\begin{array}{ll}
                 \left( \ln^2 N \right) / N^2,&i=0,\ldots,N/4-1,\\\\
                  1/ N^2,& i=N/4,\ldots,3N/4,\\\\
                 \left( \ln^2 N\right)/ N^2,& i=3N/4+1,\ldots,N,
                \end{array}   
         \right.
\end{equation*}
where $y(x_i)$ is the value of the exact solution, $\overline{y}_i$ is the value of the numerical solution  of the problem \eqref{uvod1}--\eqref{uvod3} in the mesh point $x_i,$ respectively, and $C>0$ is a constant independent of $N$ and $\varepsilon.$
 \end{theorem}
\begin{proof}\ \\
\noindent \textbf{Case $1\leqslant i\leqslant N/4-1.$} Here holds $h_{i-1}=h_i=\mathcal{O}(\varepsilon\ln N/N).$

After using $\varepsilon y''(x_k)=f(x_k,y(x_k)),$ Taylor expansion for $\cosh(\beta h_i),\,y(x_{i-1}),\,y(x_{i+1})$ and some computing we get 
\begin{align}
    Ty(x_i)=&\frac{\gamma}{\beta^2h^2_i+\mathcal{O}(\beta^4h^4_i)}\left[ \left(2+\frac{\beta^2h^2_i}{2}+\mathcal{O}(\beta^4h^4_i)\right)
     \left(y''h^2_i+\frac{y^{(iv)}(\xi^{-}_i)+y^{(iv)}(\xi^{+}_i)}{24}\right)        \right.\nonumber\\
     &-\frac{\varepsilon^2}{\gamma}\left(y''(x_{i-1})+2y''(x_i)+y''(x_{i+1}) \right)
       \left.  \left( \frac{\beta^2h^2}{2}+\mathcal{O}(\beta^4h^4_i) \right) \right].
\end{align}
Now, we apply Taylor expansion on $y''(x_{i-1}),\,y''(x_{i+1})$ and obtain 
\begin{multline}
  Ty(x_i)=\frac{\gamma}{\frac{\beta^2h^2_i}{2}+\mathcal{O}(\beta^4h^4_i)}
        \left[\left(\frac{\beta^2h^2_i}{2}+\mathcal{O}(\beta^4h^4_i)\right)
        \left(y''h^2_i+\frac{y^{(iv)}(\xi^{-}_i)+y^{(iv)}(\xi^{+}_i)}{24}h^4_i \right) \right.\\
          +\left.\frac{y^{(iv)}(\xi^{-}_i)+y^{(iv)}(\xi^{+}_i)}{12}h^4_i
          -\varepsilon^2\frac{\left(y^{(iv)}(\mu^{-}_i)+y^{(iv)}(\mu^{+}_i)\right)h^2_i}{2\gamma}
          \left( \frac{\beta^2h^2_i}{2}+\mathcal{O}(\beta^4h^4_i)\right)
            -\frac{4\varepsilon^2y''(x_i)\mathcal{O}(\beta^4h^4_i)}{\gamma} \right] ,
\end{multline}
where $\mu_i^{-},\,\xi_i^{-}\in(x_{i-1},x_i),\:\mu_i^{+},\,\xi_i^{+}\in(x_{i},x_{i+1}),$ and finally 
\begin{equation}
    \left| Ty(x_i)\right|\leqslant \left(C\ln^2 N\right) / N^2 ,\:i=1,\ldots,N/4-1.
\label{uniform1}
\end{equation}

\noindent \textbf{Case $N/4\leqslant i\leqslant N/2-1$.} 
The difference operator \eqref{operator1} can be written as
\begin{align*}
Ty(x_i)=&\frac{\gamma}{\frac{\cosh(\beta h_{i-1})-1}{\sinh(\beta h_{i-1})}+\frac{\cosh(\beta h_i)-1}{\sinh(\beta h_i)}}
   \left[\frac{\cosh(\beta h_{i-1})-1}{\sinh(\beta h_{i-1})}(y(x_{i-1})-y(x_i))-\frac{\cosh(\beta h_{i})-1}{\sinh(\beta h_{i})}(y(x_{i})-y(x_{i+1}))\right.\\
       & +\frac{2(y(x_{i-1})-y(x_i))}{\sinh(\beta h_{i-1})}-\frac{2(y(x_{i})-y(x_{i+1}))}{\sinh(\beta h_{i})}\\
       &\quad -\left.\frac{f(x_{i-1},y(x_{i-1}))+f(x_i,y(x_i))}{\gamma}\cdot \frac{\cosh(\beta h_{i-1})-1}{\sinh(\beta  h_{i-1})}
              -\frac{f(x_{i},y(x_{i}))+f(x_{i+1},y(x_{i+1}))}{\gamma}\cdot\frac{\cosh(\beta h_{i})-1}{\sinh(\beta  h_{i})}\right]
\end{align*}
Now by Lemma \ref{lema1}, Lemma \eqref{lema2}, Lemma \eqref{lema3} we have the estimate

\begin{equation}
      \left| Ty(x_i)\right|\leqslant C/ N^2 ,\:i=N/4,\ldots,N/2-1.
  \label{uniform2}
\end{equation}

\noindent \textbf{Case $i=N/2$.} It is easy to prove this case, because $h_{N/2-1}=h_{N/2},$ and the influence of the layer component $s$ is negligible, and holds
\begin{equation}
   \left|Ty(x_{N/2})\right|\leqslant C / N^2.
   \label{uniform3}
\end{equation}
Finally, according \eqref{uniform1}, \eqref{uniform2} and \eqref{uniform3}, the proof is complete.
\end{proof}

\section{Numerical experiments}
\begin{figure}[!h] 
	\includegraphics[scale=0.45]{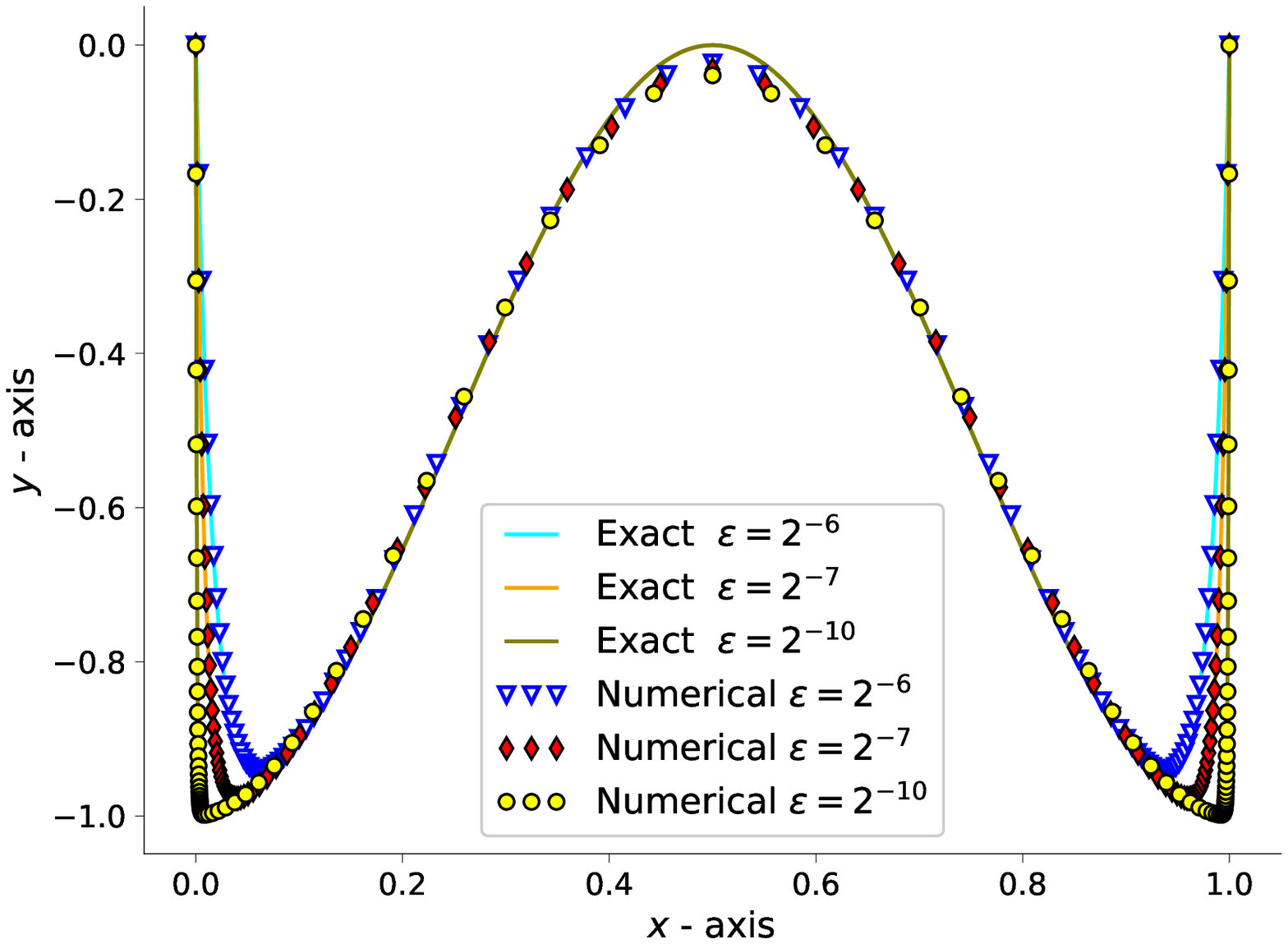}\hspace{-.5cm}	
	\includegraphics[scale=0.45]{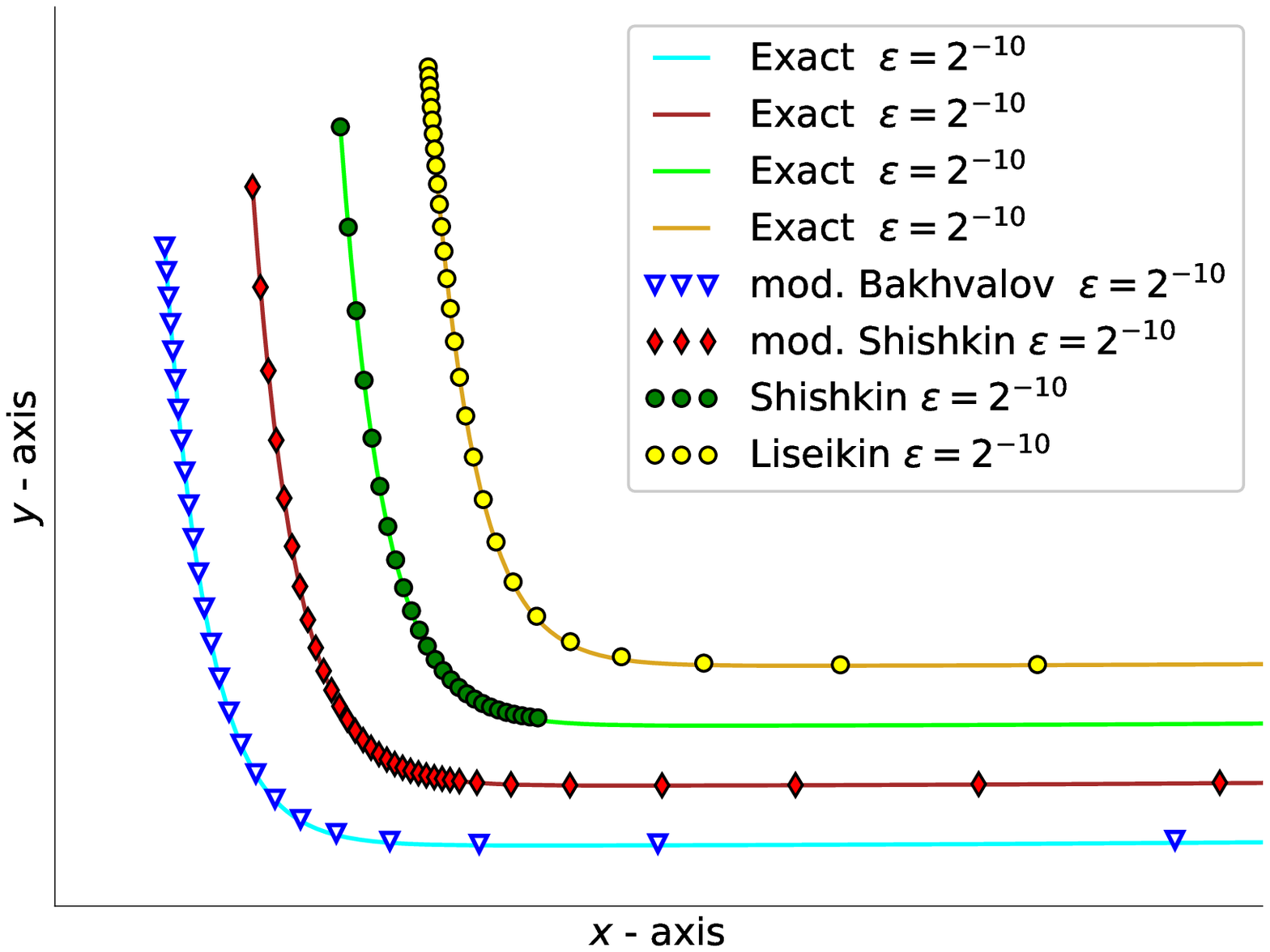}
	\caption{Exact and numerical solutions (left), layer near $x=0$ (right)--the difference scheme \eqref{konst4}}
	\label{figure1}
\end{figure}

In this section we present numerical results to confirm the accuracy of the difference scheme \eqref{konst4} using the meshes \eqref{meshShishkin1}, \eqref{meshShishkin2},  \eqref{meshVulanovic1}, and \eqref{meshLiseikin1}. 
\begin{example}
  We consider the following  boundary value problem 
      \begin{equation}
    	\varepsilon ^2y''=y+\cos^2\pi x+2\varepsilon ^2\pi^2\cos^2\pi x\quad\text{on}\quad\left(0,1\right),
    	\label{problem1}
    	\end{equation}
    	\begin{equation}
    	y(0)=y(1)=0.
    	\label{problem}
    	\end{equation}
    	The exact solution of this problem is
    	\begin{equation}
    	y(x)=\dfrac{e^{-\frac{x}{\varepsilon}}+e^{\frac{x}{\varepsilon}}}{1+e^{-\frac{1}{\varepsilon}}}-\cos^2\pi x.
    	\label{problem3}
    	\end{equation} 	
The nonlinear  system was solved using the initial condition $y_0=-0.5$ and the value of the constant $\gamma=1.$ Because the fact that the exact solution is known, we compute the error $E_N$ and the rate of convergence Ord in the usual way
\begin{equation}
   E_N=\Vert y-\overline{y}^{N}\Vert_{\infty},\: {\rm Ord}=\frac{\ln E_N-\ln E_{2N}}{\ln(2k/(k+1))},({\rm Shishkin}),\:\:{\rm Ord}=\frac{\ln E_N-\ln E_{2N}}{\ln 2},\:({\rm Bakhvalov, Liseikin})
\end{equation}
where $N=2^{k},$ $k=4,5,\ldots,12,$ and $y$ is the exact solution of the problem \eqref{uvod1}--\eqref{uvod3}, while $\overline{y}^N$ an appropriate numerical solution of \eqref{uvod1}--\eqref{uvod3}. The graphics of the numerical and exact solutions, for various values of the parameter $\varepsilon$ are on Figure  \ref{figure1} (left), while fragments of these solutions on Figure \ref{figure1} (right). The values of $E_N$ and Ord are in Tables \ref{table1}. 
\end{example}

\newpage

\begin{table*}[!h]\centering	\tiny
	\begin{tabular}{c cc cc cc cc cc cc}
		\toprule 
		
		\multicolumn{1}{c}{} 
		& \multicolumn{2}{c}{$2^{-3}$}   
		& \multicolumn{2}{c}{$2^{-5}$} 
		& \multicolumn{2}{c}{$2^{-10}$} 
		& \multicolumn{2}{c}{$2^{-20}$}  
		& \multicolumn{2}{c}{$2^{-30}$}
		& \multicolumn{2}{c}{$2^{-40}$}\\
		
		\cmidrule(r{\tabcolsep}){2-3} 
		\cmidrule(r{\tabcolsep}){4-5} 
		\cmidrule(r{\tabcolsep}){6-7}
		\cmidrule(r{\tabcolsep}){8-9}
		\cmidrule(r{\tabcolsep}){10-11}
		\cmidrule(r{\tabcolsep}){12-13}
		
		\multicolumn{1}{c}{$N$}
		& \multicolumn{1}{c}{$E_n$}
		& \multicolumn{1}{c}{Ord }
		& \multicolumn{1}{c}{$E_n$}
		& \multicolumn{1}{c}{Ord }
		& \multicolumn{1}{c}{$E_n$}
		& \multicolumn{1}{c}{Ord }
		& \multicolumn{1}{c}{$E_n$}
		& \multicolumn{1}{c}{Ord }
		& \multicolumn{1}{c}{$E_n$}
		& \multicolumn{1}{c}{Ord }\\
		\midrule
		
		$2^{4}$   &6.590e-2 &3.02 &1.469e-1 &1.83 &1.678e-1 &2.50 &1.706e-1 &2.48 &1.706e-1 &2.48 &1.706e-1 &2.48 \\ 
		$2^{5}$   &1.593e-2 &2.87 &6.209e-2 &2.74 &5.172e-2 &2.26 &5.305e-2 &2.31 &5.305e-2 &2.31 &5.305e-e &2.31 \\ 
		$2^{6}$   &3.670e-3 &2.75 &1.530e-2 &2.86 &1.627e-2 &1.98 &1.626e-2 &1.99 &1.626e-2 &1.99 &1.626e-2 &1.99 \\ 
		$2^{7}$   &8.325e-4 &2.66 &3.262e-3 &2.76 &5.576e-3 &1.98 &5.575e-3 &1.98 &5.575e-3 &1.98 &5.576e-2 &1.98 \\ 
		$2^{8}$   &1.869e-4 &2.61 &6.934e-4 &2.65 &1.836e-3 &2.00 &1.836e-3 &2.00 &1.836e-3 &2.00 &1.836e-3 &2.00 \\ 
		$2^{9}$   &4.153e-5 &2.57 &1.505e-4 &2.55 &5.820e-4 &2.00 &5.819e-4 &2.00 &5.819e-4 &2.00 &5.819e-4 &2.00 \\ 
		$2^{10}$  &9.126e-6 &2.55 &3.348e-5 &2.46 &1.797e-4 &2.00 &1.797e-4 &2.00 &1.797e-4 &2.00 &1.797e.4 &2.00 \\
		$2^{11}$  &1.981e-6 &2.53 &7.659e-6 &2.37 &5.438e-5 &2.00 &5.438e-5 &2.00 &5.438e-5 &2.00 &5.438e-5 &2.00 \\
		$2^{12}$  &4.249e-7 &-    &1.815e-6 &-    &1.618e-5 & -   &1.618e-5 & -   &1.618e-5 &-    &1.618e-5 &-  \\
		\midrule \\
		& & & & & & mesh \eqref{meshShishkin1}& & & & & &\\
		\midrule  
		
		$2^{4}$   &2.349e-1 &2.68 &3.860e-1 &1.72 &4.453e-1 &1.45 &4.339e-1 &1.42 &4.533e-1 &1.42 &4.533e-1 &1.42 \\ 
		$2^{5}$   &6.662e-2 &3.00 &1.718e-1 &2.37 &2.242e-1 &1.92 &2.317e-1 &1.89 &2.317e-1 &1.89 &2.317e-1 &1.89 \\ 
		$2^{6}$   &1.434e-2 &3.00 &5.098e-2 &2.68 &8.388e-2 &2.20 &8.802e-2 &2.18 &8.802e-2 &2.18 &8.802e-2 &2.18 \\ 
		$2^{7}$   &2.840e-3 &2.96 &1.961e-2 &2.65 &2.555e-2 &2.32 &2.717e-2 &2.30 &2.718e-2 &2.30 &2.718e-2 &2.30 \\ 
		$2^{8}$   &5.407e-4 &2.94 &2.714e-3 &2.65 &6.946e-3 &2.36 &7.486e-3 &2.33 &7.487e-3 &2.33 &7.487e-3 &2.33 \\ 
		$2^{9}$   &9.946e-5 &2.94 &6.214e-4 &2.50 &1.782e-3 &2.39 &1.955e-3 &2.32 &1.955e-3 &2.32 &1.955e-3 &2.32 \\ 
		$2^{10}$  &1.765e-5 &2.94 &1.425e-4 &2.46 &4.362e-4 &2.49 &4.987e-4 &2.30 &4.988e-4 &2.30 &4.988e-4 &2.30 \\
		$2^{11}$  &3.033e-6 &2.93 &3.263e-5 &2.43 &9.826e-5 &2.53 &1.258e-4 &2.27 &1.259e-4 &2.27 &1.259e-4 &2.27 \\
		$2^{12}$  &5.129e-7 &-    &7.446e-6 &-    &2.111e-5 &-    &3.161e-5 &-    &3.162e-5 &-    &3.162e-5 &-  \\
		\midrule \\
		& & & & & & mesh \eqref{meshShishkin2}& & & & & &\\
		\midrule
		
		$2^{4}$   &1.759e-2 &2.00 &4.723e-2 &1.94 &8.772e-2 &1.88 &1.273e-1 &1.80 &1.277e-1 &1.80 &1.277e-1 &1.80 \\ 
		$2^{5}$   &4.440e-3 &2.00 &1.226e-2 &2.02 &2.367e-2 &1.97 &3.652e-2 &1.95 &3.666e-2 &1.94 &3.666e-2 &1.94 \\ 
		$2^{6}$   &1.102e-3 &2.00 &3.010e-3 &1.97 &6.023e-3 &1.89 &9.479e-3 &1.99 &9.514e-3 &1.99 &9.515e-3 &1.99 \\ 
		$2^{7}$   &2.757e-4 &2.00 &7.696e-4 &2.00 &1.623e-3 &1.65 &2.392e-3 &2.00 &2.401e-3 &2.00 &2.401e-3 &2.00 \\ 
		$2^{8}$   &6.894e-5 &2.00 &1.923e-4 &2.00 &5.140e-4 &1.90 &5.995e-4 &2.00 &6.018e-4 &2.00 &6.019e-4 &2.00  \\ 
		$2^{9}$   &1.723e-5 &2.00 &4.807e-5 &2.00 &1.372e-4 &2.07 &1.499e-4 &2.00 &1.505e-4 &2.00 &1.505e-5 &2.00 \\ 
		$2^{10}$  &4.309e-6 &2.00 &1.201e-5 &2.00 &3.258e-5 &2.01 &3.750e-5 &2.00 &3.764e-5 &2.00 &3.764e-5 &2.00 \\
		$2^{11}$  &1.077e-6 &2.00 &3.004e-6 &2.00 &8.114e-6 &2.03 &9.375e-6 &2.00 &9.411e-6 &2.00 &9.412e-6 &2.00 \\
		$2^{12}$  &2.693e-7 & -   &7.511e-7 &-    &1.985e-6 &-    &2.343e-6 &-    &2.352e-6 &-    &2.353e-6 &- \\
		\midrule \\
		& & & & & & mesh \eqref{meshVulanovic1} & & & & & &\\

		\midrule
		
		$2^{4}$    &6.113e-2  &1.98  &1.231e-1  &1.96  &1.859e-1  &1.65 &2.014e-1 &1.63 &2.027e-1 &1.63 &2.088e-1 &1.63 \\ 
		$2^{5}$    &1.546e-2  &1.99  &3.151e-2  &2.09  &5.895e-2  &1.86 &6.497e-2 &1.85 &6.548e-2 &1.85 &6.553e-2 &1.85 \\ 
		$2^{6}$    &3.869e-3  &1.99  &7.395e-3  &2.04  &1.627e-2  &1.94 &1.795e-2 &1.94 &1.810e-2 &1.93 &1.811e-2 &1.93 \\ 
		$2^{7}$    &9.675e-4  &2.00  &1.792e-3  &2.01  &4.192e-3  &1.99 &4.678e-3 &1.97 &4.719e-3 &1.97 &4.723e-3 &1.97 \\ 
		$2^{8}$    &2.418e-4  &2.00  &4.439e-4  &2.00  &1.052e-3  &2.06 &1.191e-3 &1.98 &1.201e-3 &1.98 &1.202e-3 &1.98 \\ 
		$2^{9}$    &6.046e-5  &2.00  &1.107e-4  &2.00  &2.517e-4  &2.18 &3.004e-4 &1.99 &3.030e-4 &1.00 &3.033e-4 &1.99 \\ 
		$2^{10}$   &1.511e-5  &2.00  &2.766e-5  &2.00  &5.531e-5  &2.18 &7.542e-5 &2.00 &7.608e-5 &2.00 &7.614e-5 &2.00 \\ 
		$2^{11}$   &3.779e-6  &2.00  &6.915e-6  &2.00  &1.219e-5  &2.07 &1.889e-5 &2.00 &1.905e-5 &2.00 &1.907e-5 &2.00 \\ 
		$2^{12}$   &9.448e-7  &-     &1.728e-6  &-     &2.887e-6  &-    &4.728e-6 &-    &4.769e-6 &-    &4.773e-6 &- \\ 
	 
		\midrule \\
		& & & & & & mesh \eqref{meshLiseikin1} & & & & & &\\

	    \bottomrule
	\end{tabular}
	\caption{Results for difference scheme \eqref{konst4}}
	\label{table1}
\end{table*}

\newpage 
 
\section{Conclusion }
In this paper we give a discretization of a semilinear reaction--diffusion one--dimensional boundary--value problem. The difference scheme is constructed, the  $\varepsilon$--uniform convergence of the order 2 on modified Shishkin mesh is shown. Similar results were obtained in one of the previous papers of the first author. But in that previous paper, the approximation of the function $\psi, $ which appears in the integrals,  is very clumsily chosen.  This made the analysis unnecessarily difficult. In this paper we have chosen a simplier approximation,  which caused a simplier analysis.  Another difficulty from the mentioned previous paper we didn't overcome, here we used in our analysis the modified Shishkin mesh too, introduced by Vulanovi\'c.  It is remain a task for some future paper to replace the modified Shishkin mesh by the  Shishkin mesh.  In the numerical experiments except the modified Shishkin mesh we used the Shishkin, the modified Bakhvalov, and the Liseikin mesh. All the meshes gave the expected  results.


\printbibliography

@article{samir2010scheme,
	author = {Duvnjakovi\'{c},E. and Karasulji\'{c},S. and Oki\v{c}i\'{c},N.}, 
	title = "{Difference Scheme for Semilinear Reaction-Diffusion Problem}",
	booktitle = {14th International Research/Expert Conference Trends in the Development of Machinery 
	and Associated Technology TMT 2010, 7. Mediterranean Cruise }, 
	year = {2010},            
	pages = {793--796},
	url={http://www.tmt.unze.ba/zbornik/TMT2010/199-TMT10-138.pdf}
}

@article{samir2011scheme,
	author = { Duvnjakovi\'c, E. and Karasulji\'c, S.},
	title = {Difference Scheme for Semilinear Reaction-Diffusion Problem on a Mesh of Bakhvalov Type}, 
	journal = {Mathematica Balkanica},
	volume = {25, Fasc. 5},
	pages ={499--504},
	year = {2011},
	url={http://www.mathbalkanica.info/toc/siteCONT25-5.pdf}
}

@article{samir2015uniformlyconvergent,
	author={Duvnjakovi{\'c}, E. and Karasulji{\'c}, S. and Pa{\v s}i{\'c}, V. and Zarin, H.},
	title={A uniformly convergent difference scheme on a modified Shishkin mesh for the singularly perturbed reaction-diffusion boundary value problem},
	year={2015},
	journal={Journal of Modern Methods in Numerical Mathematics},
	volume={6},
	number={1},
	pages={28-43},
	url={http://www.m-sciences.com/index.php?journal=jmmnm&page=issue&op=view&path[]=72},
	doi={10.20454/jmmnm.2015.971}
}

@article{samir20152d,
	author={Karasulji{\'c}, S. and Duvnjakovi{\'c}, E. and  Zarin, H.},
	title={A uniformly convergent difference scheme on a Shishkin type mesh for the 2D singular perturbation boundary value problem},
	year={2015},
	note={Submitted}
}

@article{samir2015uniformly,
	author={Karasulji{\'c}, S. and Duvnjakovi{\'c}, E. and  Zarin, H.},
	title={Uniformly convergent difference scheme for a semilinear reaction-diffusion problem},
	journal={Advances in Mathematics: Scientific Journal},
	volume={4},
	number={2},
	year={2015},
	pages={139--159},
	url={http://research-publication.com/wp-content/uploads/2019/03/AMSJ-2015-N2-6.pdf}
}

@article{samir2015modifiedbah,
	author={Karasulji{\'c}, S.  and Duvnjakovi{\'c}, E. and Zarin, H.},
	title={A uniformly convergent difference scheme on a modified Bakhvalov mesh for the singular perturbation boundary value problem},
	year={2018},
	note={Submitted}
}

@article{samir2017construction,
	author={Karasulji{\' c}, S. and Duvnjakovi{\' c}, E. and Pa{\v s}i{\' c}, V. and Barakovi{\' c}, E.},
	title={Construction of a global solution for the one dimensional singularly--perturbed boundary value problem}, 
	journal={Journal of Modern Methods in Numerical Mathematics},
	volume={8},
	number={1--2},
	year={2017},
	doi={10.20454/jmmnm.2017.1275 },
	pages={52-65}    
}

@article{samir2018uniformly,
	author={Karasulji{\' c}, S. and Duvnjakovi{\' c}, E. and Memi{\' c}, E.},
	title={Uniformly Convergent Difference Scheme for a Semilinear Reaction-Diffusion Problem on a Shishkin Mesh},
	journal={Advances in Mathematics: Scientific Journal},
	volume={7},
	number={1},
	year={2018},
	pages={23--38},
	url={http://research-publication.com/wp-content/uploads/2019/03/AMSJ-2018-N1-4.pdf}
}

@article{karasuljic2019class,
	title={A class of difference schemes uniformly convergent on a modified Bakhvalov mesh},
	author={Karasulji{\'c}, S. and Zarin, H. and Duvnjakovi{\'c}, E.},
	journal={Journal of Modern Methods in Numerical Mathematics},
	volume={10},
	number={1-2},
	pages={16--35},
	year={2019},
	publisher={ModernScience Publishers},
	doi={10.20454/jmmnm.2019.1513}
}

@article{liseikin2020numerical,
	title={Numerical analysis of grid--clustering rules for problems with power of the first type boundary layers},
	author={Liseikin, V. D. and Karasulji{\'c}, S.},
	journal={Computational Technologies},
	year={2020},
	volume={25},
	number={1},
	pages={49--66},
	doi={10.25743/ICT.2020.25.1.004}
}

@article{karasuljic2020onconstruction,
	title={On construction of a global numerical solution for a semilinear singularly-perturbed reaction 
	diffusion boundary value problem},
	author={Karasulji{\'c}, S. and Ljevakovi{\'c}, H.},
	journal={Mat. bilten},
	volume={44},
	number={2},
	pages={131--148},
	year={2020},
	doi={10.37560/matbil2020131k}
}

@incollection{Liseikin_2021,
	doi = {10.1007/978-3-030-76798-3_14},
	url = {https://doi.org/10.1007/2F978-3-030-76798-3_14},	
	year = {2021},	
	publisher = {Springer International Publishing},	
	pages = {227--240},	
	author = {V. D. Liseikin and S. Karasuljic and A. V. Mukhortov and V. I. Paasonen},	
	title = {On a Comprehensive Grid for Solving Problems Having Exponential or Power-of-First-Type Layers},	
	booktitle = {Lecture Notes in Computational Science and Engineering}
}

@inproceedings{samir2010tmt, 
	author = {Duvnjakovi\'{c},E. and Karasulji\'{c},S. and Oki\v{c}i\'{c},N.}, 
	title = {Difference Scheme for Semilinear Reaction-Diffusion Problem},
	booktitle = {14th International Research/Expert Conference Trends in the Development of Machinery 
	and Associated Technology TMT 2010, 7. Mediterranean Cruise }, 
	year = {2010},            
	pages = {793--796},
	url={ http://www.tmt.unze.ba/zbornik/TMT2010/199-TMT10-138.pdf}
}

@inproceedings{samir2011uniformnly,
	author = {Duvnjakovi\'{c},E. and Karasulji\'{c},S.},
	Title={Uniformly Convergente Difference Scheme for Semilinear Reaction-Diffusion Problem },
	booktitle = {Conference on Appllied and Scietific Computing, Trogir, Croatia },
	pages = {25},
	year = {2011 },
	url={http://applmath11.math.hr/abs_book.pd}
}

@inproceedings{samir2011skoplje,
	author = {Duvnjakovi\'{c},E.  and Karasulji\'{c},S.},
	Title={Uniformly Convergente Difference Scheme for Semilinear Reaction-Diffusion Problem  },
	booktitle = {SEE Doctoral Year Evaluation Workshop, Skopje, Macedonia},
	year = {2011 },
	url={http://users.sch.gr/alexiouth/index.php/workshops/evaluation-ws}
}

@inproceedings{samir2012uniformnly,
	author = {Karasulji\'{c},S. and Duvnjakovi\' {c},E.},
	Title={Difference Scheme for Semilinear Reaction-Diffusion Problem on Layer--adapted Mesh},
	booktitle = {The Seventh Bosnian-Herzegovinian Mathematical Conference, Sarajevo, BiH},
	year = {2012 },
	url={http://www.anubih.ba/Journals/vol.8,no-2,y12/16matskup7.pdf}
}

@inproceedings{samir2012class,
	author = {Duvnjakovi\'{c},E.  and Karasulji\'{c},S.},
	Title={Class of Difference Scheme for Semilinear Reaction-Diffusion Problem on Shishkin Mesh},
	booktitle = {MASSEE International Congress on Mathematics - MICOM 2012, Sarajevo, Bosnia and Herzegovina},
	year = {2012 }
}

@inproceedings{samir2013collocation,
	author = {Duvnjakovi\'{c},E. and Karasulji\'{c},S.},
	Title={Collocation Spline Methods for Semilinear Reaction-Diffusion Problem on Shishkin Mesh},
	booktitle = {IECMSA-2013, Second International Eurasian Conference on Mathematical Sciences and Applications, Sarajevo, Bosnia and Herzegovina},
	year = {2013 }
}

@inproceedings{samir2015construction,
	author = {Karasulji{\'{c}},S. },
	Title={Construction of the Difference Scheme for Semilinear Reaction-Diffusion  Problem on a Bakhvalov Type Mesh},
	booktitle = {The Eighth Bosnian-Herzegovinian Mathematical Conference, Sarajevo, BiH},
	year = {2015 },
	% url={http://www.anubih.ba/Journals/vol.11,no-2,y15/14Resume.pdf}
}

@inproceedings{liseikin2019rules,
	title={On Rules for Grid Clustering in the Zones of Boundary and Interior Layers},
	author={Liseikin, V.D. and Kudryavtsev, A. N. and Paasonen, V. I. and Karasulji{\'c}, S. and Mukhortov, A. V.},
	booktitle={Mathematics and its Applications. International Conference  in honor of the 90th birthday of Sergei K. Godunov},
	year={2019},
	organization={Novosibirsk, Russia}
}

@inproceedings{liseikin2020comprehensiveproceedings,
	title={On a comprehensive grid for solving problems having exponential or power--of--first--type layers},
	author={Liseikin, Vladimir and Karasulji{\'c}, Samir and Mukhortov, Aleksandr and Paasonen, Viktor},
	booktitle={NUMGRID 2020,  Moscow, Russia},
	year={2020},
	organization={Dorodnicyn Computing Center FRC CSC RAS}
}

@book{samir2021matematika1,
	title={Matematika 1},
	author={Karasulji\'c, S. and Halilovi\'c, S.},
	year={2021},
	publisher={Off-set Tuzla},
	isbn={978-9958-31-482-7}	   
}

@book{liseikin2021numerical,
	title={Numerical Grids and High-Order Schemes for Problems with Boundary and Interior Layers},
	author={Liseikin, V.D. and Karasulji\'c, Samir and Paasonen, V.I.},
	year={2021},
	publisher={Novosibirsk State University},
	isbn={978-5-4437-1274-1}
}

@article{boglaev1984approximate,
       title = "Approximate solution of a non-linear boundary value problem with a small parameter for the highest-order differential ",
       journal={ Zh. Vychisl. Mat. i Mat. Fiz.},
       volume = "24",
       number = "11",
       pages = "1649 -- 1656",
       year = "1984",
       doi = {10.1016/0041-5553(84)90005-3},
       url = "http://www.mathnet.ru/links/3d5582cba89c929035e5d7e541089694/zvmmf4282.pdf",
       author = {Boglaev, I. P.},
       note={In Russian}
       }

@article{herceg1990,
       year={1989},
       journal={Numer. Math.},
       volume={56},
       number={7},
       doi={10.1007/BF01405196},
       title={Uniform fourth order difference scheme for a singular perturbation problem},
       url={http://dx.doi.org/10.1007/BF01405196},
       publisher={Springer-Verlag},
       author={Herceg, D.},
       pages={675--693},
       url = {http://eudml.org/doc/133421}
       }

@article{herceg1991,
        author={Herceg,D. and Surla,K.},
        title={Solving a nonlocal singularly perturbed problem by spline in tension},
        journal={Review of research Faculty of Sciences-University of Novi Sad},
        volume={Vol.21, No.2},
        pages={119-132},
        year={1991},
        url={http://www.dmi.uns.ac.rs/nsjom/Papers/21_2/NSJOM_21_2_119_132.pdf}
        }

@article{herceg2003,
        title={On numerical solution of semilinear singular perturbation problems by using the Hermite scheme on a new Bakhvalov-type mesh},
        author={Herceg, D. and Miloradovi\'{c}, M.},
        journal={Novi Sad J. Math},
        volume={33},
        number={1},
        pages={145--162},
        year={2003},
        url={http://www.dmi.uns.ac.rs/nsjom/Papers/33_1/nsjom_33_1_145_162.pdf}
        }

@article{herceg2003a,
        title={On a fourth-order finite difference method for nonlinear two-point boundary value problems},
        author={Herceg, D. and Herceg, {Dj}.},
        journal={Novi Sad J. Math},
        volume={33},
        number={2},
        pages={173--180},
        year={2003},
        url={http://www.dmi.uns.ac.rs/nsjom/Papers/33_2/nsjom_33_2_173_180.pdf}
        }

@article{kopteva2001,
  title={Uniform second-order pointwise convergence of a central difference approximation for a quasilinear convection-diffusion problem},
  author={Kopteva, N. and Lin{\ss}, T.},
  journal={J. Comput. Appl. Math.},
  volume={137},
  number={2},
  pages={257--267},
  year={2001},
  publisher={Elsevier},
  doi={10.1016/S0377-0427(01)00353-3}
}

@article{kopteva2001robust,
  title={A robust adaptive method for a quasi-linear one-dimensional convection-diffusion problem},
  author={Kopteva, N. and Stynes, M.},
  journal={SIAM Journal on Numerical Analysis},
  volume={39},
  number={4},
  pages={1446--1467},
  year={2001},
  publisher={SIAM},
  doi={10.1137/S003614290138471X}
}

@article{kopteva2004,
  title={Numerical analysis of a singularly perturbed nonlinear reaction--diffusion problem with multiple solutions},
  author={Kopteva, N. and Stynes, M.},
  journal={Appl. Numer. Math.},
  volume={51},
  number={2},
  pages={273--288},
  year={2004},
  publisher={Elsevier},
  doi={10.1016/j.apnum.2004.07.001}
}

@article{kopteva2007maximum,
  title={Maximum norm error analysis of a 2d singularly perturbed semilinear reaction-diffusion problem},
  author={Kopteva, N.},
  journal={Mathematics of Computation},
  volume={76},
  number={258},
  pages={631--646},
  year={2007},
  doi={10.1090/S0025-5718-06-01938-7 },
  url={http://www.ams.org/journals/mcom/2007-76-258/S0025-5718-06-01938-7},
  publisher={Washington, DC: National Academy of Sciences-National Research Council,[1960?-}
}

@article{kopteva2009,
  title={A robust overlapping Schwarz method for a singularly perturbed semilinear reaction-diffusion problem with multiple solutions},
  author={Kopteva, N. and Pickett, M. and Purtill, H.},
  journal={Int. J. Numer. Anal. Model},
  volume={6},
  pages={680--695},
  year={2009},
  url={http://www.global-sci.org/ijnam/readabs.php?vol=6&no=4&doc=680&year=2009&ppage=695}
}

@article{liseikin2019compact,
        author={Liseikin, V.D. and Poaasonen},
        title={Compact Difference Schemes and Layer Resolving Grids for Numerical Modeling of Problems with Boundary and Interior Layers},
        journal={Numer. Analys. Appl.},
        volume={12},
        pages={37--50},
        year={2019},
        doi={10.1134/S199542391901004X}        
        }

@article{linss2000uniform,
  title={Uniform pointwise convergence on Shishkin-type meshes for quasi-linear convection-diffusion problems},
  author={Lin{\ss}, T. and Roos, H.G. and Vulanovi{\'c}, R.},
  journal={SIAM J. Numer. Anal.},
  volume={38},
  number={3},
  pages={897--912},
  year={2000},
  publisher={SIAM},
  doi={10.1137/S0036142999355957}
       }

@article{linss2012approximation,
  title={Approximation of singularly perturbed reaction-diffusion problems by quadratic $C^1$-splines},
  author={Lin{\ss}, T. and Radojev, G. and Zarin, H.},
  journal={Numerical Algorithms},
  volume={61},
  number={1},
  pages={35--55},
  year={2012},
  publisher={Springer US},
  doi={10.1007/s11075-011-9529-7}
         }

@article{lorenz1982stability,
         author={Lorenz, J.},
         title={Stability and monotonicity properties of stiff quasilinear boundary problems},
         journal={Zb.rad. Prir. Mat. Fak. Univ. Novom Sadu, Ser. Mat.},
         volume={ 12},
         year={1982},
         pages={ 151--176},
         note={ MR 85e:34046},
         url={http://www.emis.de/journals/NSJOM/Papers/12/NSJOM_12_151_175.pdf}
         }

@article{shishkin1988grid,
  title={Grid approximation of singularly perturbed parabolic equations with internal layers},
  author={Shishkin, G. I.},
  journal={Sov. J. Numer. Anal. M.Russian Journal of Numerical Analysis and Mathematical Modelling},
  volume={3},
  number={5},
  pages={393--408},
  year={1988},
  language={ In Russian},
  doi={ 10.1515/rnam.1988.3.5.393} 
}

@article{stynes1996,
         author = {Sun, G. and Stynes, M.},
         journal = {Math. Comput.},
         number ={ 215},
         pages = {1085--1109},
         title = {A uniformly convergent method for a singularly perturbed semilinear reaction-diffusion problem with multiple solutions},
         volume = {65},
         year = {1996},
         url={http://www.jstor.org/stable/2153793}
        }

@article{stynes2006numerical,
  title={Numerical analysis of singularly perturbed nonlinear reaction-diffusion problems with multiple solutions},
  author={Stynes, M. and Kopteva, N.},
  journal={Computers and Mathematics with Applications},
  volume={51},
  number={5},
  pages={857--864},
  year={2006},
  publisher={Elsevier},
  doi={10.1016/j.camwa.2006.03.013}
}

@article{surla2003,
         author={ Surla, K. and Uzelac, Z.},
         title={On Stability of Spline Difference Scheme for Reaction-Diffusion Time-Dependent Singularly Perturbed Problem},
         journal={Novi Sad J. Math.},
         volume={33},
         number={2},
         year={2003},
         pages={89-94},
         url={http://www.dmi.uns.ac.rs/nsjom/Papers/33_2/nsjom_33_2_089_094.pdf}
         }

@article{vulanovic1983numerical,
         author={Vulanovi{\' c}, R.},
         title={On a Numerical Solution of a Type of Singularly Perturbed Boundary Value Problem by Using a Special Discretization Mesh},
         journal={Novi Sad J. Math.},
         volume={13},
         year={1983},
         pages={187-201},
         url={http://www.dmi.uns.ac.rs/nsjom/Papers/13/NSJOM_13_187_201.pdf}
         }

@article{vulanovic1983,
        author={Vulanovi\'c, R.},
        title={On a numerical solution of a type of singularly perturbed boundary value problem by using a special discretization mesh},
        journal={Univ. u Novom Sadu Zb. Rad, Prirod-Mat. Fak. Ser. Mat},
        volume={13},
        year={1983},
        pages={187--201},
        url={http://www.dmi.pmf.uns.ac.rs/nsjom/Papers/13/NSJOM_13_187_201.pdf}
       }

@article{vulanovic1989,
        author={Vulanovi\'c, R.},
        title={Mesh generation methods for numerical solution of quasilinear singular perturbation problems},
        journal={Univ. u Novom Sadu Zb. Rad, Prirod-Mat. Fak. Ser. Mat},
        volume={19},
        number={2},
        year={1989},
         pages={171--193},
         url={http://www.emis.ams.org/journals/NSJOM/Papers/19_2/NSJOM_19_2_171_193.pdf}
       }

@article{vulanovic1991second,
  title={A second order numerical method for non-linear singular perturbation problems without turning points},
  author={Vulanovi{\'c}, R.},
  journal={USSR Comp. Math. Math+.},
  volume={31},
  number={4},
  pages={522--532},
  year={1991}
}

@article{vulanovic1993,
        author={Vulanovi\'c, R.},
        title={On numerical solution of semilinear singular perturbation problems by using the Hermite scheme},
        journal={Univ. u Novom Sadu Zb. Rad, Prirod-Mat. Fak. Ser. Mat},
        volume={23},
        number={2},
        year={1993},
        pages={363--379},
        url={http://www.dmi.uns.ac.rs/nsjom/Papers/23_2/NSJOM_23_2_363_379.pdf}
       }

@article{vulanovic2001higher,
 author = {Vulanovi\'{c}, R.},
 title = {A Higher-order Scheme for Quasilinear Boundary Value Problems with Two Small Parameters},
 journal = {Computing},
 volume = {67},
 number = {4},
 pages = {287--303},
 url = {http://dx.doi.org/10.1007/s006070170002},
 doi = {10.1007/s006070170002},
 acmid = {566399},
 publisher = {Springer-Verlag New York, Inc.},
 year={2001}
}

@article{vulanovic2004 ,
        author={Vulanovi\'c, R.},
        title={An almost sixth-order finite-difference method for semilinear singular perturbation problems},
        journal={Computational methods in applied  mathematics},
        volume={4},
        number={3},
        year={2004},
        pages={368--383},
        doi={10.2478/cmam-2004-0020}
       }

@book{liseikin2018grid,
      author={Liseikin, V. D.},
      title={Grid Generation for Problems with Boundary and Interior Layers},
      publisher={Novosibirsk State University, Novosibirsk}, 
      year={2018},
      isbn={978-5-4437-0822-5}
      }

@book{ortega2000,
        author={Ortega,J. M. and Rheinboldt,W. C.},
        title={Iterative Solution of Nonlinear Equations in Several Variables},
        publisher={SIAM, Philadelphia, USA},
        year={2000},
        isbn={9780898714616}
         }
\end{document}